\def\timestamp{%
Time-stamp: <no-categoricity.tex: Monday 16-01-2006 at 14:30:10 (cet)>}
\def\stripname Time-stamp: <#1 #2>{#2}
\edef\filedate{\expandafter\stripname\timestamp}
\documentclass
               {amsart}

\usepackage{amsrefs}

\theoremstyle{plain}
\newtheorem{theorem}{Theorem}[section]
\newtheorem{proposition}[theorem]{Proposition}
\newtheorem{lemma}[theorem]{Lemma}

\theoremstyle{remark}
\newtheorem{remark}[theorem]{Remark}

\newcommand{\B}{\mathcal{B}}
\newcommand{\C}{\mathcal{C}}
\newcommand{\D}{\mathcal{D}}
\newcommand{\U}{\mathcal{U}}
\newcommand{\V}{\mathcal{V}}

\newcommand\Hyper[1]{2^{#1}}
\newcommand\abs[1]{\mathopen|#1\mathclose|}

\DeclareMathSymbol\restr{2}{AMSa}{"16}
\DeclareMathSymbol\N{0}{AMSb}{`N}
\DeclareMathSymbol\Q{0}{AMSb}{`Q}

\begin{document}

\title{There is no categorical metric continuum}

\author{Klaas Pieter Hart}
\address{Faculty of Electrical Engineering, Mathematics, and
         Computer Science\\
         TU Delft\\
         Postbus 5031\\
         2600~GA {} Delft\\
         the Netherlands}
\email{K.P.Hart@TUDelft.NL}
\urladdr{http://fa.its.tudelft.nl/\~{}hart}

\begin{abstract} 
We show there is no categorical metric continuum.
This means that for every metric continuum $X$ there is another metric 
continuum $Y$ such that $X$ and $Y$ have (countable) elementarily
equivalent bases but $X$ and $Y$ are not homeomorphic.
As an application we show that the chainability of the pseudoarc is
not a first-order property of its lattice of closed sets.
\end{abstract}

\date{\filedate}


\maketitle

\section*{Introduction}

Many properties of compact Hausdorff spaces can, naturally, be phrased in terms
of their families of closed sets.
For a fair number of these one can find even first-order formulas in the
language of lattices that characterize them, see, e.g., 
\cite{vanderSteeg2003}.

In \cite{BanakhBankstonRainesRuitenburg} and 
\cite{HartvanderSteeg} it was shown
that chainability is not a first-order property.
In an earlier version of the former paper the question was raised whether
there is any chainable continuum for which its chainability is expressible
in first-order terms.
The authors offered the pseudoarc as a candidate.

If the pseudoarc were `first-order chainable' then it would at once become
a categorical continuum.
This is so because the pseudoarc is the only continuum that is both chainable
and hereditarily indecomposable.
Therefore any continuum with a lattice-base for its closed sets
that is elementarily equivalent to some lattice base for the closed sets
of the pseudoarc would itself be the pseudoarc.

In this note we show that no categorical continuum exists and hence,
indirectly, that the pseudoarc is not first-order chainable.

\section{Preliminaries}

\subsection{Categoricity}

Categoricity is a model-theoretic notion; we refer to 
\cite{HodgesShorterModelTheory}*{Section~6.3} 
for a complete treatment of the countable case,
which is the case that we shall need;
we refer to \cite{HodgesShorterModelTheory} for other model-theoretic
notions as well.
A countable structure $S$ (group, lattice, ordered set) is categorical if
every other countable structure that satisfies the same first-order sentences
as~$S$ is actually isomorphic to~$S$.
A prime example is the set $\Q$ of rational numbers; it is, up to 
isomorphism, the only countable linearly ordered set that is densely ordered
and without end~points, see~\cite{Cantor1895}*{\S~9}.
Structures that satisfy the same first-order sentences are usually said
to be \emph{elementarily equivalent}.

We extend these notions to cover compact Hausdorff spaces:
we call two such spaces elementarily equivalent if they have bases for
the closed sets that are elementarily equivalent as lattices.
A compact metric space is categorical if every compact metric space that
is elementarily equivalent to it is homeomorphic to it. 

As an example we mention the Cantor set: if $X$ is compact metric and if it 
has a countable base that is elementarily equivalent to some base for the
Cantor set then one readily shows that 1)~$X$ has no isolated points and 
2)~$X$ is zero-dimensional; therefore $X$~is homeomorphic to the Cantor set.

\subsection{Ultrapowers}

We use ultrapowers to find structures that are elementarily equivalent
to a given structure but, in a well-defined way, much richer.
If $L$~is a lattice and $u$~is an ultrafilter on the set~$\N$ of natural 
numbers then the ultrapower of~$L$ by the ultrafilter~$u$ is the 
power~$L^\N$ modulo the equivalence relation $x\equiv_u y$, defined
by $x\equiv_u y$ if{}f $\{n:x(n)=y(n)\}\in u$.
We denote this quotient structure by~$L_u$.
See {Section~8.5} of~\cite{HodgesShorterModelTheory} for more information
on ultraproducts and for the definition of `richness' alluded to above.

\subsection{Creating surjections}

The following lemma is used to construct continuous surjections.

\begin{lemma}[\cite{DowHart}*{Theorem 1.2}]\label{lemma.how.to.map.onto}
Let $X$ and $Y$ be compact Hausdorff spaces and let $\C$~be a base for the
closed subsets of\/~$Y$ that is closed under finite unions and finite
intersections.
Then $Y$~is a continuous image of\/~$X$ if and only if there is a map
$\phi:\C\to\Hyper{X}$ such that
\begin{enumerate}
\item $\phi(\emptyset)=\emptyset$ and\label{cond.i}
      if $F\neq\emptyset$ then $\phi(F)\neq\emptyset$;
\item if $F\cup G=Y$ then $\phi(F)\cup\phi(G)=X$; and\label{cond.ii}
\item if $F_1\cap\cdots\cap F_n=\emptyset$\label{cond.iii}
      then $\phi(F_1)\cap\cdots\cap\phi(F_n)=\emptyset$.\qed
\end{enumerate}
\end{lemma}

\subsection{$K_0$-functions}

Consider a metric space $X$, with metric $d$, and a closed subspace~$A$.
Define a map
$\kappa:\Hyper{A}\to\Hyper{X}$ by
$$
\kappa(F)=\bigl\{x\in X: d(x,F)\le d(x,A\setminus F)\bigr\}.
$$
In \cite{Kuratowski66}*{\S\,21\,XI} it is shown that for all closed
sets $F$ and~$G$ in~$A$ we have
\begin{itemize}
\item $\kappa(F)\cap A=F$;
\item $\kappa(F\cup G)=\kappa(F)\cup\kappa(G)$; and
\item $\kappa(A)=X$ and $\kappa(\emptyset)=\emptyset$ ---
      by the fact that $d(x,\emptyset)=\infty$ for all~$x$.
\end{itemize}
Following \cite{vanDouwen-thesis} we call such a function a 
\emph{$K_0$-function}.

\subsection{Chainability}

A continuum is \emph{chainable} if every finite open cover has a finite chain
refinement, that is, an indexed refinement $\{V_i:i<n\}$ with the property 
that $V_i\cap V_j\neq\emptyset$ if and only if $\abs{i-j}\le 1$.
The condition that $\V$~is a chain refinement of~$\U$
can be expressed by a (rather long) first-order formula.
The condition that $\U$~has a chain refinement is,
a priori, not first-order as one does not know beforehand how large
the refinement is going to be.
One gets a formula of the form 
$(\exists\V)\bigl(\bigvee_n\phi_n(\U,\V)\bigr)$, where $\phi_n$
expresses that $\V$~is an $n$-element chain refinement of~$\U$
--- this is an $L_{\omega_1,\omega}$-formula: each $\phi_n$~is first-order
but the disjunction is infinite.
Chainability proper is then defined by infinitely many such formulas:
one for each possible cardinality of~$\U$.

The authors of \cite{BanakhBankstonRainesRuitenburg} identified one way of 
defining first-order chainability: make sure the disjunction becomes
finite.
This would mean, in words: for every natural number~$m$ there
is a natural number~$n$ such that every open cover of size~$m$
has an open chain refinement of size~$n$ or less.

The negation of this, namely that there is a natural number~$m$ such that
for every~$n$ there is an open cover for which every chain~refinement
has at least $n$~members, was called \emph{elastically chainable}
in~\cite{BanakhBankstonRainesRuitenburg}.
However, Theorem~4.1 of~\cite{BanakhBankstonRainesRuitenburg} implies that
this is not a new property: is $X$~is a connected and normal space then
for every~$n$ it has a three-element open cover with no chain refinement
of size~$n$ or less.

Another result announced in~\cite{BanakhBankstonRainesRuitenburg} is that
the infinite number of formulas given above can be reduced to one:
a continuum is chainable if{}f every four-element open cover has 
a chain refinement.

\section{The main lemma}

\begin{lemma}\label{lemma.main}
Let $X$ and $Y$ be metric continua and let $\B$ and $\C$ be a countable
lattice bases for their respective families of closed sets.
Let $u$ be any free ultrafilter on~$\omega$.
There is a map~$\phi$ from~$\C$ to the ultrapower~$\B_u$ that satisfies
the conditions in Lemma~\ref{lemma.how.to.map.onto}.
\end{lemma}

\begin{proof}
We consider $Y$ embedded in the Hilbert cube $Q$ and we let
$\kappa:2^Y\to2^Q$ be a $K_0$-function.
Furthermore, fix a continuous surjection $f:X\to[0,1]$.

Enumerate $\C$ as $\langle C_n:n\in\omega\rangle$ and put
$E=\{e\subseteq\omega:\bigcap_{i\in e}C_i=\emptyset\}$.
Observe that $Y\cap\bigcap_{i\in e}\kappa(C_i)=\emptyset$ whenever $e\in E$.

Fix $n<\omega$ and take a positive number $\epsilon_n$ less than~$2^{-n}$ 
and all distances between~$Y$ and $\bigcap_{i\in e}\kappa(C_i)$ for
those $e\in E$ that are subsets of~$n$.
Take a continuous map $g_n:[0,1]\to Q$ such that the image is a subset of
$B(Y,\epsilon_n)$ and such that it meets every ball $B(y,\epsilon_n)$
with $y\in Y$ (here we use that $Y$ is a continuum: it has arbitrarily
small arcwise connected neighbourhoods).

For $i<n$ let $D^n_i$ be the preimage 
$f^\gets\bigl[g_n^\gets[\kappa(C_i)]\bigr]$.
Because $\kappa$~is a $K_0$-function we know that
$D^n_i\cup D^n_j=X$ whenever $C_i\cup C_j=Y$.
Also, by the choice of~$\epsilon_n$, we know that 
$\bigcap_{i\in e}D^n_i=\emptyset$ whenever $e\in E$ and $e\subseteq n$.
Now expand the sets $D^n_i$ to get members $B^n_i$ of~$\B$, 
retaining the property that 
$\bigcap_{i\in e}B^n_i=\emptyset$ whenever $e\in E$ and $e\subseteq n$.

The definition of~$\phi$ is now straightforward:
define $\phi(C_i)$ to be the $\equiv_u$-equivalence class of
$\langle B^n_i:n>i\rangle$.
Note that $\phi$ has the required properties even when we take the reduced
power modulo the co-finite filter.
\end{proof}

\section{The main result}

The following proposition is the key to the main result.

\begin{proposition}\label{prop.preimage}
Let $X$ and $Z$ be two metric continua.
There is a third metric continuum~$Y$ such that 
\begin{enumerate}
\item $Z$~is a continuous image of~$Y$; and
\item $Y$ and $X$ have elementarily equivalent bases for the closed sets.
\end{enumerate}
\end{proposition}

\begin{proof}
Take countable bases $\B$ and $\D$ respectively for the closed sets
of~$X$ and $Z$.
Fix a free ultrafilter~$u$ on~$\omega$ and apply 
Lemma~\ref{lemma.main} to find a map $\phi:\D\to \B_u$ as in
Lemma~\ref{lemma.how.to.map.onto}.
Next apply the L\"owenheim-Skolem theorem to obtain a countable elementary
substructure~$\C$ of~$\B_u$ that contains $\phi[\D]$.
We let $Y$ be the Wallman space of the lattice~$\C$.
Then $Y$~is as required: the lattice $\C$ is elementarily equivalent
to~$\B_u$ and hence to~$\B$ itself.
The map~$\phi$ enables us, via Lemma~\ref{lemma.how.to.map.onto},
to map~$Y$ onto~$Z$.
\end{proof}

\subsection{The proof}

It is now straightforward to prove the main assertion of this note.
In \cite{Waraszkiewicz34} Wa\-rasz\-kie\-wicz constructed a family of continua
such that no single metric continuum maps onto all of them.
Let $X$ be any metric continuum and fix a continuum~$Z$ from that family
that is not a continuous image of~$X$.
Apply Proposition~\ref{prop.preimage} to find a metric continuum~$Y$
that does map onto~$Z$ and yet has a base for the closed sets that
is elementarily equivalent to a base for the closed sets of~$X$.
Clearly $X$ and $Y$ are not homeomorphic.

\begin{remark}
The referee observed that the main result remains valid if 
`compact metric' is replaced by `compact and of weight less 
than~$2^{\aleph_0}$'.
Indeed, the proof in~\cite{Waraszkiewicz34}
establishes that if $X$~is any continuum that maps onto all continua
in the family constructed there then the space $C(X,\mathbb{R}^2)$ 
(with the uniform metric) has a discrete subspace of 
cardinality~$2^{\aleph_0}$, in fact there is a constant~$a$ such that
if $f$ and $g$ map $X$ onto different members of the family then
their uniform distance is at least~$a$.
This implies that $X$ can be mapped onto at most $w(X)$ many members
of the family.

No essential modifications are needed.
One should observe that in Lemma~\ref{lemma.main} 
and Proposition~\ref{prop.preimage} the continuum~$X$ need not be metric
and in the latter proposition one can take $Y$~to be of the same weight 
as~$X$.
\end{remark}

\subsection{The pseudoarc}

In an earlier version of~\cite{BanakhBankstonRainesRuitenburg} it was asked
whether the pseudoarc is inelastically chainable.
If it were it would show that the pseudoarc is categorical.

The results of this paper imply that this corollary does not hold
and hence that the pseudoarc \emph{is} elastically chainable.
This argument simply shows that a natural number~$m$ as in the definition
exists, it does not provide a definite value.

Of course this particular corollary has been superseded by results
from~\cite{BanakhBankstonRainesRuitenburg}; as noted above
\emph{every} connected normal space is elastically chainable in the sense 
that for every natural number~$N$ there is a three-element open cover 
that cannot be refined by a chain-cover with fewer than $N$ elements.

\end{document}